# A Correlation Analysis and Visualization of Climate Change using Post-Disaster Heterogeneous Datasets

**Sukeerthi Mandyam**
Computer Science and Engineering Department
Sri Sivasubramaniya Nadar College of Engineering
Email: sukeerthi18175@cse.ssn.edu.in

**Shanmuga Priya M G**
Computer Science and Engineering Department
Sri Sivasubramaniya Nadar College of Engineering
Email: shanmugapriya18140@cse.ssn.edu.in

**Shalini S**
Computer Science and Engineering Department
Sri Sivasubramaniya Nadar College of Engineering
Email: shalini18139@cse.ssn.edu.in

**Kavitha S**
Computer Science and Engineering Department
Sri Sivasubramaniya Nadar College of Engineering
Email: kavithas@ssn.edu.in

Abstract: There are numerous geo-climatic and human factors that contribute to the occurrence of natural disasters in the real-world scenario. Besides the study of causes and preconditions of such calamities, post-disaster analysis is essential for the efficient management of the disaster situation. This process needs timely and accurate data in light of the increasing frequency and severity of climate change-related extreme weather events. The analysis of disaster data involves the challenging task of integrating multiple heterogeneous sources, data ingestion and visualization. This paper aims at providing a three-dimensional analytical view of disaster data as time-series charts and a statistical model to evaluate the correlation between the occurrence of disasters, climate change and the corresponding economic damages (as a percentage of GDP). The created corpus consists of two unified heterogeneous datasets. The first one forms the basis for exploratory data analysis involving visualization charts for disaster-related factors. The second dataset involves temperature anomaly information which allows for climate change analysis. The integration of the above datasets provides scope for correlation analysis using three different variations of coefficients. Therefore, statistical methodologies are leveraged to play important role in managing disasters, from preparation to recovery and reporting. The graphical representations provide insights on regional trends that follow, related to factors such as the proportion of each type of disaster in the various losses incurred. Therefore, obtaining reliable information about the population, the economy and climate are crucial both for risk management and preparedness and for responding to disasters. The research puts forth a detailed statistical methodology with a spatial dimension to study the impacts on the economy and infrastructure in the aftermath of a disaster thereby ensuring specialized assistance. The inference of the analysis confirmed a positive correlation between climate change and occurrence of natural disasters. Therefore, statistical evidence of an important phenomenon like climate change affecting natural disasters brings awareness among the population in the society to be more environmentally responsible.

## 1. Introduction

Natural disasters have long been an area of great interest in the research community since they repeatedly cause damage and economic loss in various fields. The current scope of research cannot predict the occurrence and severe impacts of a disaster with acceptable accuracy. The issue of limited resources emerges because this unpredictability makes it difficult to allocate adequate resources in advance. Although the effects of natural disasters on humans are extensive, economic development is hindered to a large extent years after their occurrence.

Environmentalists suggested that these erratic patterns are probably due to climate change and it is high time we transit quickly to a preparedness-centric approach rather than continuing to be in the relief-centric mode. The inferences of relevant patterns between climate change and disaster occurrence and corresponding economic damages helps the population and the government to become more aware of the causes of the disasters that can be controlled through proper action. The analysis also provides appropriate metrics for efficient planning post the disaster occurrence to ensure faster relief and emergency response through the allocation of required resources.

The creation of disaster data corpus from heterogeneous datasets unifies disparate databases and provides a uniform conceptual schema to build a common interface. By integrating closely-related datasets using common attributes, a consolidated analysis is possible by incorporating multiple factors. For exploratory data analysis of these disaster factors, the main dataset was outsourced from EM-DAT (Emergency Events Database), which is a global database on natural and technological disasters, containing essential core data on the occurrence and effects of more than 21,000 disasters in the world, from 1900 to the present [1]. By extracting the useful metrics for analysis from raw datasets and combining them with data from open-source platforms like Kaggle offering the earth's temperature over time, an integrated corpus was created to meet the requirements for the analysis methodology [2].

A considerable number of inferences can be fetched from the dataset by defining certain metrics as attributes for the analysis. For instance, time-series charts for visualizing the disaster scenario can provide a whole new dimensional view of the data. To gain an insightful perspective about the raw data, it needs to be combined with relatable effects such as economic losses. Statistical formulas can be used to focus on an alarming aspect related to disasters such as climate change and find out the actual consequential impact of global warming. Therefore, a correlation analysis is formulated between the average temperature of the earth and the frequency of natural disasters and the corresponding economic indicators of a region.

The remaining sections of the paper are structured as follows: Section 2 discusses the literature survey involved in disaster management, followed by Section 3 which describes the methods adopted in analysis. Section 4 is subdivided as experimental setup and performance analysis elucidating the results which are followed by concluding remarks in Section 5.

## 2. Literature Survey

The research on existing literature initially focuses on identifying the causes of disasters and gathering divergent perceptions specifically about climate change. We observe the trends in disaster damages related to life and economic losses and how effective early warning systems have proven to be in emergency-response situations. Furthermore, we explore empirical assessments that emphasize the way climate change translates to economic damages and the importance of including vulnerability and socio-economic factors for analysis. The data integration methods and the latest open-source technologies needed to build a consistent data corpus and visualizations were studied. Lastly, we examine suitable statistical approaches that could be used for analysis such as the correlation between factors and a comparative model for inferring results across various regions and types of disasters.

UNECE has made recommendations to help countries harness official statistics to improve disaster response [3]. They have narrowed down two stark trends that are becoming ever-more present in today's world: rapid urbanization, and a growing number and intensity of hazardous events related to climate change. The combination of these two forces results in increased exposure to risk for people, the economy and the environment. A recent study conducted by Gonzalo et al. found that people have different perceptions of climate change and its effects on disaster occurrences that differ regionally and socio-economically [4]. The studies related to disasters are performed to estimate damage sizes by using various medium variables such as previous damage data, and economic indicators by countries. UNDRR (UN Office for Disaster Risk Reduction) published a report to mark the International Day for Disaster Risk Reduction on October 13, 2020, and this confirms the extent to which extreme weather and climate change have dominated the disaster landscape in recent times [5].

The World Meteorological Organization informed that weather-related disasters have increased causing more damage but fewer deaths. The number of disasters has increased by a factor of five over the 50-year period, driven by climate change, more extreme weather and improved reporting [6]. Fortunately, due to improved early warning mechanisms, the number of deaths has significantly gone down. From 1970 to 2020, weather, climate and water hazards accounted for 50% of all disasters, 45% of reported deaths and 74% of reported economic losses. Since there is considerable scope for improved analysis of disaster aspects, we aim to provide a more specific analysis of the dependency between these factors using relevant visualization and statistical methods.

The study by Coronese et al. presented evidence for an increase in the economic damages due to natural disasters [7]. They have analyzed event-level data using regressions to capture patterns in the economic damage distribution. This pattern is strongest in temperate regions, suggesting that the prevalence of devastating natural disasters has broadened beyond tropical regions and that adaptation measures have had some mitigating effects on damages. Geiger et al. expressed serious concerns regarding the analysis and the interpretation of the results [8]. They suggested that the capital stock metric should be used as the dependent variable instead of GDP which requires exclusive access. Their interpretation was that climate change was more attributed to the intensity of the disaster events rather than economic damages. This particular research focused only on the shortcomings of the previous study and does not provide a methodology to overcome them. Therefore, we aim to apply a correlative approach to infer the dependency between the three factors with open-source data.

Donner et al. proposed that the socio-economic diversity and demographics of a region play an important role in assessing the disaster risk and vulnerability of a region [9]. SAMHSA published their findings on how disasters affect people of low socioeconomic status (SES) [10]. They are comparatively less prepared for disasters than others and it relates to the fact that people of low SES cannot always afford more expensive preparedness actions. Hence, we have incorporated an analysis among countries with varying sustainable development index (SDI) where lower SDI denotes countries with low standards of living.

James proposed a method for the data integration of heterogeneous datasets [11]. As more analytical functions increase, the need for interacting with heterogeneous data increases. Being able to discover data from other sources, rapidly integrate that data with one's own, and perform simple analyses, often by visualizations, leads to insights that are considered important assets. Traditional enterprise approaches such as creating data warehouses with carefully analyzed data organization and regular updates and curation are used for organizing structured data. This research proposes three major steps to achieve the objective, namely data discovery, integration and exploration but does not emphasize the suitable data warehouse to be chosen for the respective use case. By taking the feasible inputs from this paper, we have adopted a cloud data warehouse called Snowflake for our data integration needs after procuring disaster datasets from multiple sources to create our own data corpus.

Snowflake has out-of-the-box features like separation of storage and computing, data sharing and third-party business intelligence tools support for advanced analytics [12]. For modelling the charts, Apache Superset is a recently developed tool that is fast, lightweight, intuitive, and loaded with options that make it easy for users of all skill sets to explore and visualize their data even at a petabyte-scale [13]. The representation of the correlation factor as a matrix uses matplotlib and seaborn libraries as they handle the combination of data frames efficiently to provide a multi-dimensional analysis [14].

Weixiao et al. proposed a study describing the characteristics and losses related to natural disasters in China between the years 1985 and 2014 [15]. The dataset was outsourced from EM-DAT but the region was confined to China and the types of disasters analyzed were predominantly floods and storms. The damages caused due to disasters focused only on life losses and economic crisis. We improvise this analysis in multiple facets by involving abundant disaster types, losses and regions and also including a correlation analysis with climate change. The visualization charts are also presented in an interactive and easily interpretable manner using current technologies.

Sandra et al. conducted research to analyze the impact of climate change on natural disasters [16]. The study also outlines the consequences of natural disasters involving a global distribution along with future projections about weather-related events. Some of the inferences drawn highlighted that though disasters are inevitable and will always be a part of nature, their frequency and intensity have been altered due to climate change, thus making regions more susceptible to multi-faceted losses.

Suryanto proposed a method for analyzing the correlation between vulnerable areas and residents' risk perception of earthquakes in a particular province called Bantul [17]. This research focused on mapping the vulnerability of the population that potentially is affected by the negative impact of earthquakes and showing the correlation between risk perception, social variables, and economic variables to the vulnerability of the region. The perception index value was calculated based on the average of the perception of the earthquake impact, the confidence level in the earthquake-resistant housing and control capabilities. Further, the correlation method was used to determine the relationship between the individual's perception of the disaster risk on the level of the vulnerability of the region and represented the same using a colour scheme. Therefore, refining this correlative model allows us to apply a similar kind of correlative analysis from a worldwide perspective for multiple disasters using a different set of factors such as economic impacts and climate change.

Thomas et al. proposed a method to observe the relationship between climate change and increase in natural disasters in the Asia and Pacific regions [18]. The growing frequency of intense natural disasters is being attributed to three key disaster risk factors: increased population exposure, susceptibility and climate-related risks. There is a substantial relationship between the growth in natural disasters and population exposure, as indicated by population densities formulated using regression analysis. Along with the scientific link between greenhouse gases and climate change, the data in this research point to a link between rising natural disasters and human-caused greenhouse gas emissions in the atmosphere.

Michael proposed a method for analyzing the relationship between climate change, natural disasters and migration to find out any empirical evidence relating to the factors [19]. It was observed that there was a positive correlation between the frequency of occurrence of natural disasters and corresponding emigration rates of the region. The study also adds another perspective that emigration is not common in under-developed nations. The situation of migration differs widely attributing to the reason, be it a sudden occurrence of natural disaster or gradual climate change. There exists scope for applying such analysis on a broader perspective relating to economic damages caused by climate change.

Rodkin published a collection of research topics in the domain of Statistical Analysis of Natural Disasters and Related Losses [20]. They drew theoretical conclusions by applying described techniques to different natural disasters. A comparative study of losses from earthquakes, floods, tornadoes and hurricanes was presented. The losses of different types such as fatalities, number of affected people, and overall economic losses are analyzed region-wise namely in Japan and USA. The Kolmogorov test was used as the statistical tool for testing hypotheses on distribution under study. The study compares two regions with different economic strata and geographical parameters, hence the interpretations are not specialized and correlated by using socio-economic factors and historical risks of the region. However, by making use of the statistical interpretation from this publication, exploratory and relational analysis concerning multiple regions and factors are applied and results are inferred.

Therefore, by reviewing the methodologies adopted and types of inferences made in these studies, we have created a corpus including two unified heterogeneous datasets with necessary disaster and climate change factors. We have further explored a scope of three-dimensional analysis using charts and correlation evaluation by considering important aspects of disasters that can be generalized across all regions.

## 3. Methods

This section describes the methodologies adopted to develop the heterogeneous datasets required for creating visualization charts and performing correlation analysis of disaster factors with climate change.

*3.1. Description of the dataset*

A wide range of factors contributes to the occurrence of disasters as well as the impact of their consequences. A collection of five datasets was outsourced, namely, frequency of occurrence of natural disasters, details of affected population due to disasters, economic damages caused by disasters, news coverage of disasters and Earth's global temperature data.

The datasets were built into two main corpuses using Snowflake [EM-DAT (2020), Berkeley Earth (2017)]. The attribute selection from the raw data was performed by considering the usefulness of the metric, redundancy and consistency of the data. Some of the disaster-related factors are not included due to presence of a higher proportion of null values in the time frame chosen for analysis, which when considered would lead to inaccurate dependencies. The main objective of developing two heterogeneous datasets is to combine their respective analysis using common metrics while also maintaining a separate corpus for closely related factors.

The first dataset – ***Exploratory Analysis for Disaster-Related Factors (EADRF)***- consists of the distribution of disaster data by region and type of disaster with indicators for damages represented numerically as deaths, affected, injured, homeless, loss in GDP etc. The types of disasters included in the dataset are floods, earthquakes, droughts, landslides, volcanic activity, wildfires, extreme temperature and weather. The second dataset- ***Climate Change Analysis using Temperature Anomaly (CCATA)*** includes consolidated data of global temperature anomaly for correlation analysis with disaster occurrences and economic damage.

*3.1. Data Visualization*

Raw historical data on disaster types and effects cannot be utilized to make interpretations. However, insightful visualization tools of the raw disaster data can provide meaningful inferences by accounting for multiple facets. Besides using some basic statistical representations such as bar graphs and pie charts, we have included some use-case specific charts such as time-series charts and area charts. For world map visualizations, corresponding ISO codes were embedded in the dataset.

A time-series chart illustrates data points at successive intervals of time. We have chosen the appropriate interval to be yearly, to represent annual trends in disaster factors. A sunburst chart is used to conveniently visualize a hierarchical dataset without filtering, through a series of concentric rings, where each ring is segmented proportionally to represent its constituent ratio or percentage. A dual-line chart is used to represent multiple measures with two different axis ranges. An area chart 'stacks' trends on top of each other to illustrate how the percentage of certain aspects of the dataset is changing over time. Therefore, when combined with a table calculation that represents the share of the total for each dimension member in the visualization, stacked area charts are an effective way to evaluate trends in distributions.

The EADRF dataset is mainly used for developing charts related to disaster factors such as the frequency of occurrence, type and region-wise losses, economic damages and crisis-related analysis and news coverage of disasters. Whereas, the visualizations consisting of temperature anomaly and correlation analysis are created using the integrated EADRF and CCATA datasets by choosing common attributes.

*3.2. Correlation Analysis*

Disaster-related information consists of data from multidisciplinary sources involving numerous factors in the real-time scenario. To analyze this dynamic and large quantitative data, correlation analysis is a suitable method that estimates the linear association of various disaster factors. The inferences obtained from this analysis provide scope for forecasting disaster events based on the relationship with correlated entities like climate change and temperature anomaly. The correlation model also adds a dimension of statistically proving common notions of causes and effects related to disasters.

The correlation coefficient represents the strength of the relationship between the relative dynamics of two variables as a statistical measure. The values range between -1.0 and 1.0. Correlation coefficient values less than +0.8 or greater than -0.8 are not considered significant. The most commonly used correlation coefficient is the Pearson correlation (r). This measures the strength and direction of the linear relationship between two variables, therefore there is a need for a linear relationship between the pair of variables. High linearity with a straight-line graphical representation in the relationship of variables indicates high correlation between them. Equation (1) formulates the Pearson linear correlation coefficient.

$$r = \frac{\sum_{i=1}^{M} (X_i - \overline{X})(Y_i - \overline{Y})}{\sqrt{\sum_{i=1}^{M} (X_i - \overline{X})^2 (Y_i - \overline{Y})^2}} \qquad (1)$$

Spearman correlation coefficient is a statistical method of measuring the non-parametric correlation that indicates the dependency between two samples. It is applied in the cases of the lack of normally distributed nature within the variable set, therefore not requiring a linear relationship. It does not assume the distribution of data, measures variables with ranking or the size rank-order standard, and calculates the size of rank-order correlation coefficients between two variables. Spearman's method is more sensitive to error and discrepancies in the data. Equation (2) formulates the Spearman rank-order correlation coefficient.

$$r_s = 1 - \frac{6 \sum d_i^2}{N(N^2 - 1)} \qquad (2)$$

The Kendall correlation coefficient measures the non-parametric correlation that indicates the dependency between two samples using information converted with data standards or ranking standards. The interpretation of Kendall's Tau coefficient in terms of probabilities of observing agreeable and non-agreeable pairs is very direct. Equation (3) formulates the Kendall rank-order correlation coefficient.

$$\tau = \frac{N_c - N_d}{\frac{1}{2}N(N - 1)} \qquad (3)$$

The integration of data has allowed for the creation of an efficient data corpus and ingesting it through the Snowflake data warehouse enables cloud-storage of the data. The connection with Apache Superset provides an interface to develop visualization charts and consolidate the results as a dashboard for presentation. To implement the correlation analysis module, machine learning libraries are used to observe suitable patterns in data. The appropriate correlation function is chosen for the data samples and the coefficients are tabulated which is fine-tuned for better representation as a correlation heatmap matrix.

## 4. Results and Performance Analysis

*4.1. Experimental Setup*

The implementation of this research requires a minimum of 8 GB RAM to support applications like Jupyter Notebook. Apache Superset is a visualization tool that requires an Ubuntu environment for installation. To ensure minimal consumption of cloud data warehouse credits, it is advised to use an extra small warehouse size, which is ideal in the long run.

*4.2. Performance Analysis*

The two main corpora developed were further used in creating the dashboard consisting of charts for various analyses [EADRF, CCATA]. In the EADRF dataset, we have classified the disasters by type and region as shown in Fig. 1 and Fig. 2. The dataset includes information categorized across 175 countries and 8 types of disasters. This allows for a clear demarcation to formulate the disaster losses specific to each disaster and country and analyze them individually.

| ENTITY | CODE | YEAR | DEATHS | DEATH_RATE | PERCENTAGE_SHARE_DEATHS.. | INTERNALLY_DISPLACED_POPULATION |
|---|---|---|---|---|---|---|
| India | IND | 2008-01-01 | 1734.947159 | 0.143342031 | 0.019412573 | 6662000 |
| India | IND | 2009-01-01 | 1650.188873 | 0.134166941 | 0.018182388 | 5304000 |
| India | IND | 2010-01-01 | 990.4368701 | 0.079265138 | 0.010694988 | 1411000 |
| India | IND | 2011-01-01 | 824.7629692 | 0.064971695 | 0.008770145 | 1503000 |
| India | IND | 2012-01-01 | 334.1529185 | 0.025904882 | 0.003536581 | 9110000 |
| India | IND | 2013-01-01 | 6556.513259 | 0.500271762 | 0.068948811 | 2145000 |
| India | IND | 2014-01-01 | 866.6232915 | 0.065157421 | 0.009064289 | 3428000 |
| India | IND | 2015-01-01 | 1121.625116 | 0.083212952 | 0.011607983 | 3655000 |
| India | IND | 2016-01-01 | 807.587968 | 0.059180022 | 0.00821994 | 2400000 |

Fig. 1. Data classified by 175 regions consisting of 7 columns and 6469 rows

| ENTITY ↓ | YEAR | DEATHS | AFFECTED | HOMELESS | ↓ INJURED |
|---|---|---|---|---|---|
| Flood | 1982-01-01 | 4648 | 36917037 | 372410 | 25292 |
| Flood | 2015-01-01 | 3495 | 27293725 | 165658 | 23218 |
| Flood | 1994-01-01 | 6771 | 122546263 | 7214123 | 22785 |
| Flood | 2004-01-01 | 6982 | 116517821 | 457117 | 15877 |
| Flood | 1965-01-01 | 1401 | 4410813 | 88214 | 15245 |
| Flood | 1989-01-01 | 4716 | 103446717 | 867158 | 11312 |
| Flood | 2010-01-01 | 8356 | 188113195 | 670720 | 10383 |
| Flood | 2016-01-01 | 4720 | 76086888 | 2253027 | 8936 |
| Flood | 2012-01-01 | 3544 | 63730563 | 222538 | 8918 |

Fig. 2. Data classified by 8 types consisting of 6 columns and 757 rows

Some of the visualizations that constitute the dashboard built using EADRF dataset are shown as sample cases. However, there is extensive scope for creating refined visualizations filtered by either region or type for specific analytical purposes. By default, the charts represent global data unless the country is specified. The life loss incurred due to disaster is visualized in a region-wise manner. It is followed by a proportional analysis of the affected population and the extent to which economic status can play a role in disaster survival. The colour-schemed map in Fig. 3 describes the region-wise death count due to disasters. The sunburst chart in Fig. 4 represents the hierarchical proportion of deaths from the total affected population due to each type of disaster. The multi-layered time-series chart in Fig. 5 shows the trend of death rates (number of deaths per 100,000 population) of countries grouped by the SDI metric. The legend allows for a selective visualization or explicit comparison between any desired entities on the chart.

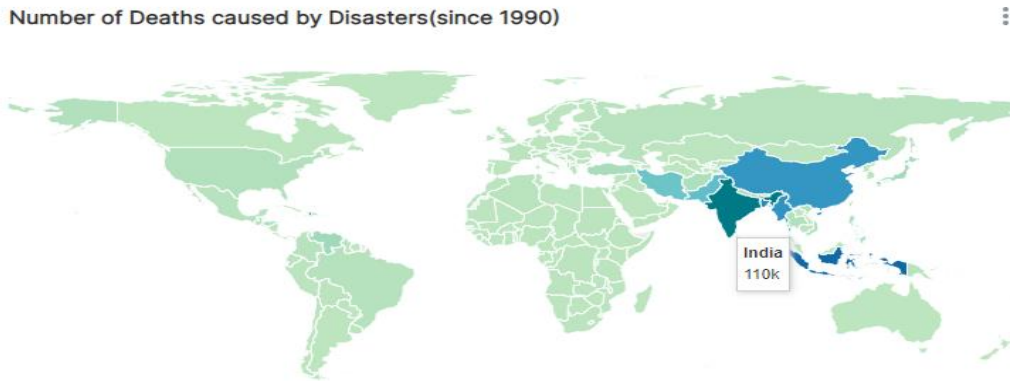

Fig. 3. A colour-schemed world-map representing the number of deaths due to disasters

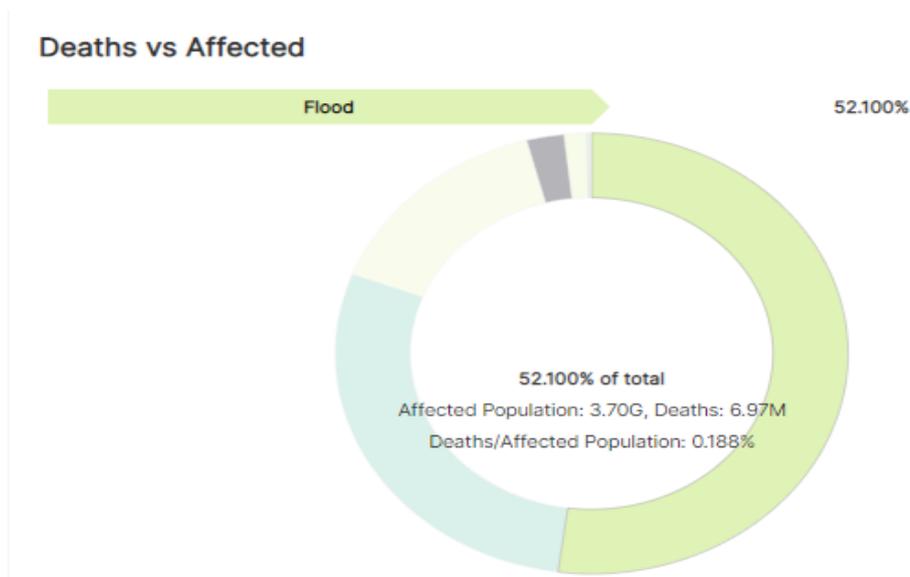

Fig. 4. A sunburst chart representing the ratio of deaths to affected population for each type of disaster

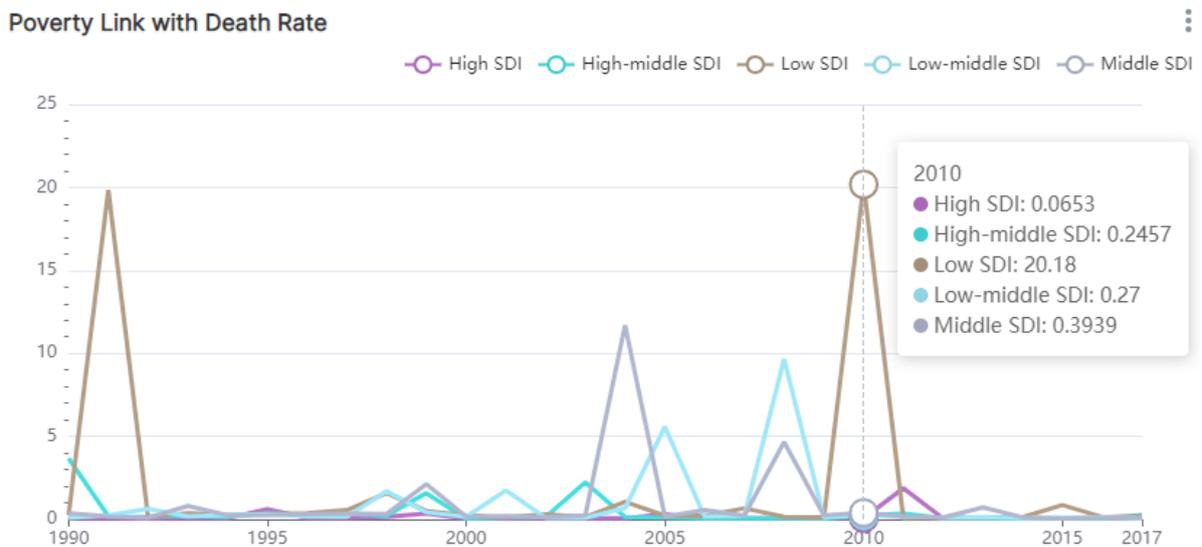

Fig. 5. A time-series chart of death rates grouped by countries based on Sustainable Development Index

Since natural disasters have devastating impacts in terms of loss of life, property and economic damage, we visualize global direct disaster losses as a share of GDP in Fig. 6.

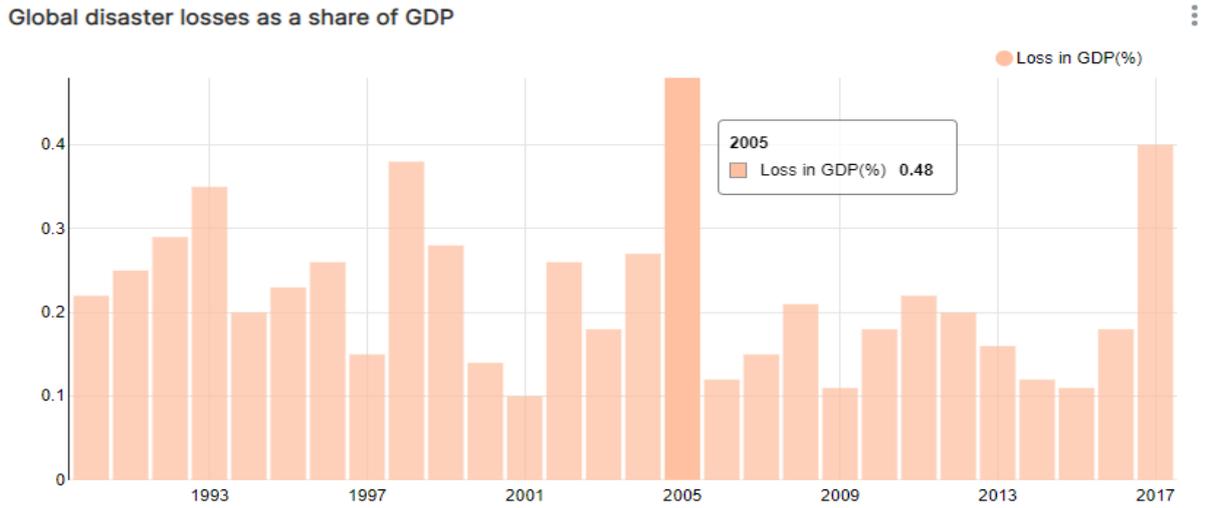

Fig. 6. A bar chart of percentage loss in GDP caused by disasters

The visualization of news coverage of disaster types shows the perception people have about its influence. The bar chart in Fig. 7 shows the percentage of news coverage for each type of disaster to rank them according to newsworthiness.

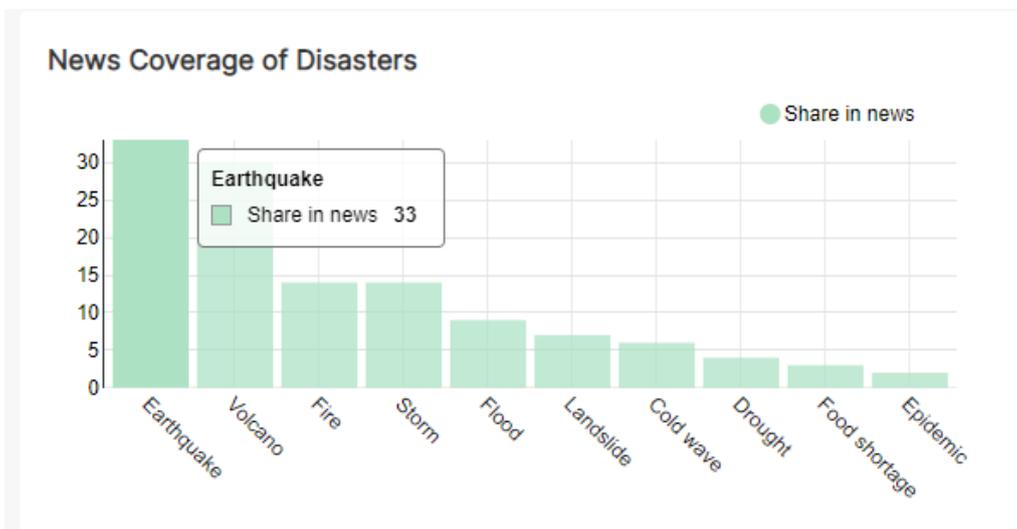

Fig. 7. A bar chart of proportion of each disaster in News Coverage

Some of the common consequences of disaster occurrence are migration, homelessness and injuries. The time-series chart in Fig. 8 shows the number of internally displaced populations across the neighbouring countries of India as well as the USA and UK. The area chart in Fig. 9 shows the proportion of affected, homeless and injured people due to disasters. The respective definitions are given as-
Injured: People suffering from physical injuries, trauma or an illness requiring immediate medical assistance.
Homeless: People whose house is destroyed or heavily damaged and therefore need shelter.
Affected: People requiring immediate assistance i.e., requiring basic survival needs such as food, water, shelter

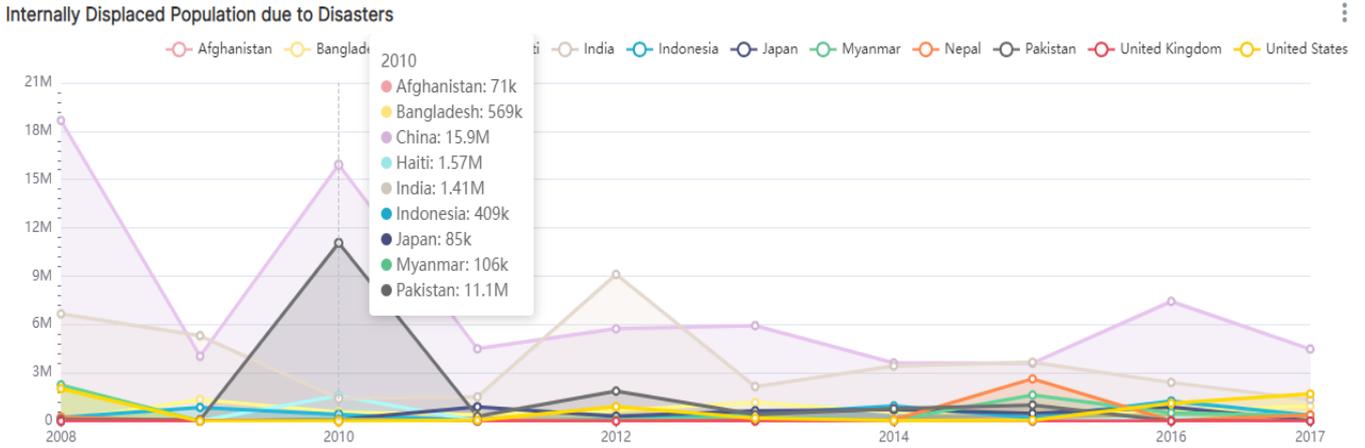

Fig. 8. A time-series chart of internally displaced population across neighbouring countries of India

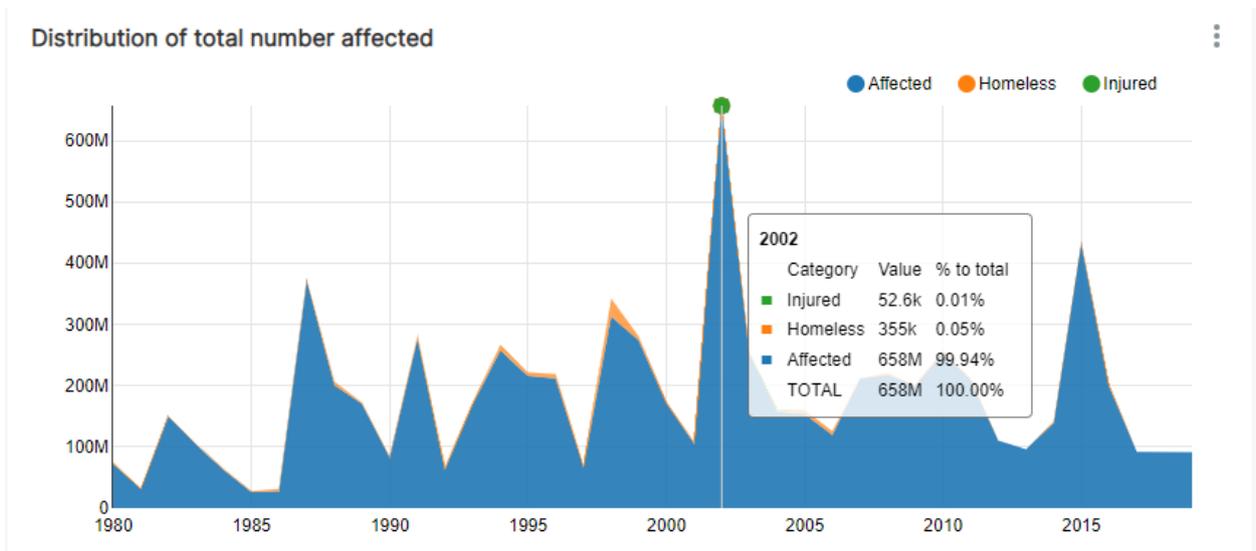

Fig. 9. An area chart consisting of in-depth analysis of total affected population

Concerning the correlation analysis module, initially, we used graphical representations plotted using raw data to find suitable patterns and scope to proceed further. The area chart shown in Fig. 10 built using EADRF, stacks up each disaster entity type indicating their frequency of occurrence over the years. It is observed that the proportion of occurrence of floods has increased significantly when compared to other disasters. A general consensus of increased risk is also seen due to rise in natural calamities as time progresses.

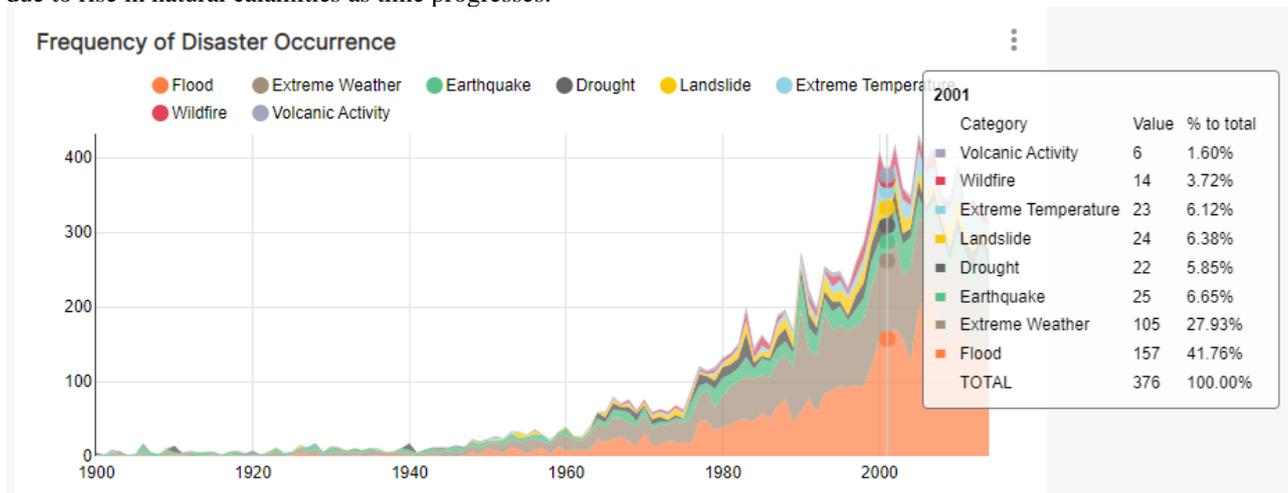

Fig. 10. An area chart of global occurrences of natural disasters between 1900 and 2015

Temperature Anomaly can be defined as the amount of deviation from a certain frame of reference, in this scenario the long-term average global temperature. It is an efficient way to observe key trends in climate change because gathering absolute temperatures across the world is a difficult task. A positive anomaly value indicates that current temperature is higher than reference level whereas a negative value indicates a cooler temperature with respect to baseline. A scientific fact that puts this into perspective is that the Earth's temperature has been rising by 0.14 Fahrenheit per decade and the warming rate since 1980 has been recorded at an alarming rate of 0.35 Fahrenheit. To estimate the feasibility of performing correlation analysis, it is necessary to plot the trend in variables to observe dependencies and their increase with time. The dual line chart shown in Fig. 11, uses data integrated from EADRF and CCATA using year as the common attribute. The trend indicates that there is a positive dependency of correlation between frequency of natural disaster occurrence and temperature deviation which is more predominant after the 1980s.

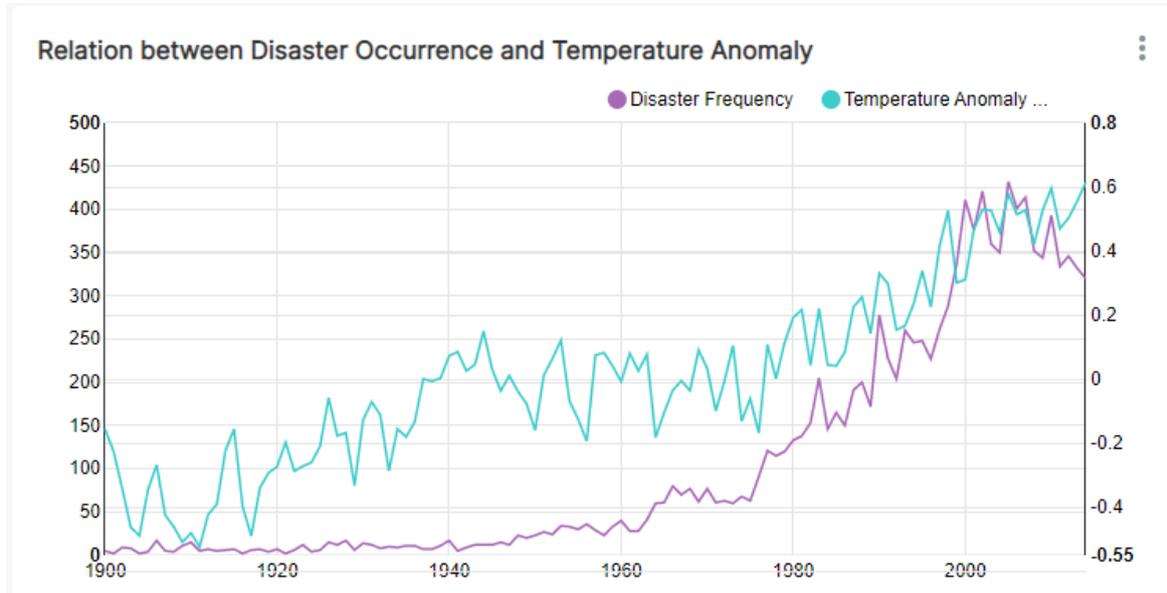

Fig. 11. A dual line chart representing all natural disaster occurrences vs temperature anomaly between 1900 and 2015

Economic damage varies greatly depending on a variety of factors such as the type of disaster, the location of the disaster, the intensity of the disaster, disaster management activities, and so on. However, one crucial tendency that emerges from the stacked area chart in Fig. 12 created using EADRF is that the total economic loss caused by all disasters has been steadily growing over time. The main factor is an increase in disaster incidence (which is linked to global warming), albeit not all disasters cause equal economic harm, for instance floods account to the most damage.

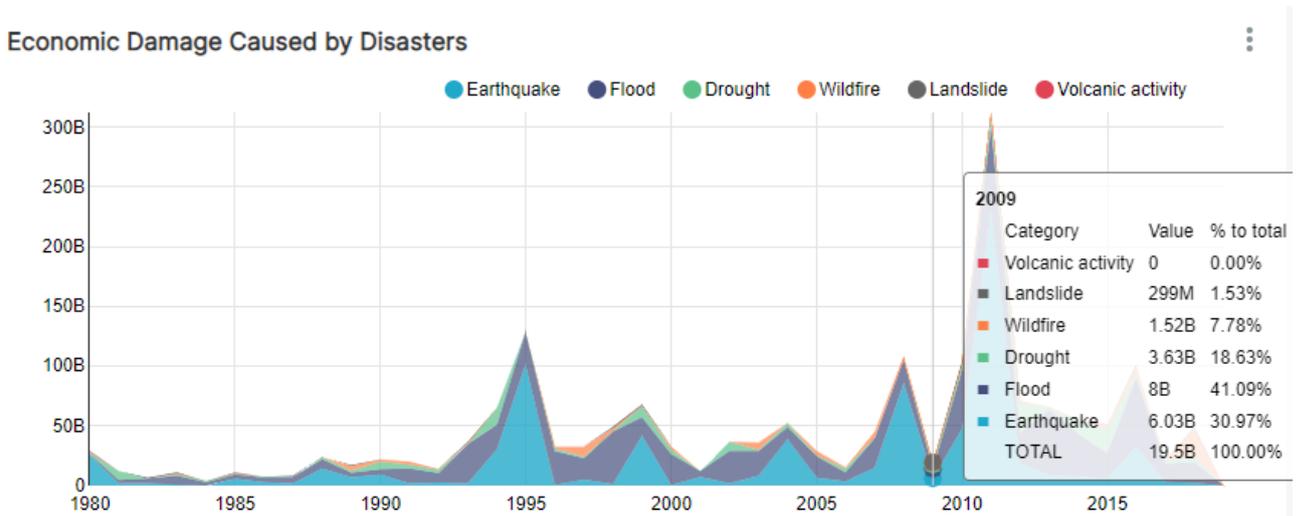

Fig. 12. An area chart of distribution of each disaster in economic damages caused

By combining the two heterogeneous datasets EADRF and CCATA, a correlation table is created to enumerate the correlation between each entity. The tables shown in Fig. 13 and Fig. 14 consist of 10 rows and 10 columns, represent the Pearson correlation coefficient between temperature anomaly with disaster occurrence and economic damage respectively. Similarly, the tabulated results and matrices are obtained for Kendall and Spearman's methodologies. The tabulated results have been limited to the head of the data frame while choosing a colour-schemed correlation matrix representation for the complete result of the analysis.

|  | Temperature Anomaly | All natural disasters | Drought | Earthquake | Extreme temperature | Extreme weather | Flood | Landslide | Volcanic activity | Wildfire |
|---|---|---|---|---|---|---|---|---|---|---|
| Temperature Anomaly | 1.000000 | 0.865128 | 0.750086 | 0.796451 | 0.738191 | 0.848536 | 0.841740 | 0.772180 | 0.669120 | 0.683839 |
| All natural disasters | 0.865128 | 1.000000 | 0.854789 | 0.903392 | 0.859863 | 0.969517 | 0.978128 | 0.894296 | 0.788590 | 0.837557 |
| Drought | 0.750086 | 0.854789 | 1.000000 | 0.815286 | 0.684009 | 0.829133 | 0.785536 | 0.790511 | 0.710884 | 0.788781 |
| Earthquake | 0.796451 | 0.903392 | 0.815286 | 1.000000 | 0.721105 | 0.917687 | 0.824128 | 0.778272 | 0.708440 | 0.778318 |
| Extreme temperature | 0.738191 | 0.859863 | 0.684009 | 0.721105 | 1.000000 | 0.776985 | 0.870008 | 0.751340 | 0.556775 | 0.694524 |

Fig. 13. Pearson correlation table between temperature anomaly and disaster occurrence

|  | Temperature Anomaly | All natural disasters | Drought | Earthquake | Extreme temperature | Extreme weather | Flood | Landslide | Volcanic activity | Wildfire |
|---|---|---|---|---|---|---|---|---|---|---|
| Temperature Anomaly | 1.000000 | 0.647406 | 0.515849 | 0.346882 | 0.294110 | 0.588786 | 0.678628 | 0.353444 | 0.215112 | 0.557517 |
| All natural disasters | 0.647406 | 1.000000 | 0.498355 | 0.823475 | 0.373682 | 0.791300 | 0.832297 | 0.244076 | 0.103445 | 0.504414 |
| Drought | 0.515849 | 0.498355 | 1.000000 | 0.270173 | 0.026217 | 0.415827 | 0.466098 | 0.126578 | 0.247509 | 0.242899 |
| Earthquake | 0.346882 | 0.823475 | 0.270173 | 1.000000 | 0.291325 | 0.364361 | 0.624735 | 0.083710 | 0.087658 | 0.261585 |
| Extreme temperature | 0.294110 | 0.373682 | 0.026217 | 0.291325 | 1.000000 | 0.273596 | 0.251023 | 0.061126 | 0.022204 | 0.352994 |

Fig. 14. Pearson correlation table between temperature anomaly and economic damages

Pearson correlation heatmaps for visualizing the dependency of temperature anomaly with disaster occurrence and economic damage are shown in Fig. 15 and Fig. 16 respectively. The heatmap obtained supports our initial assumption that disasters are indeed a consequence of climate change as indicated by the strong positive correlation coefficient value. Since economic damages constitute broader real-word issues rather than being confined to climate change, they exhibit relatively lower correlation with temperature anomaly.

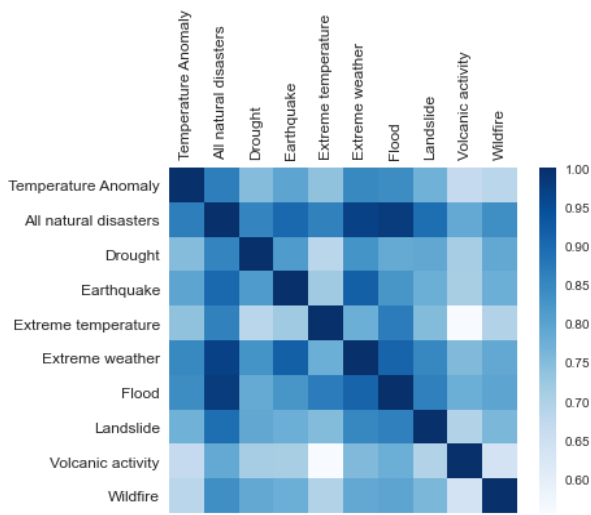

Fig. 15. Pearson Correlation Heatmap for dependencies between disaster occurrence and temperature anomaly.

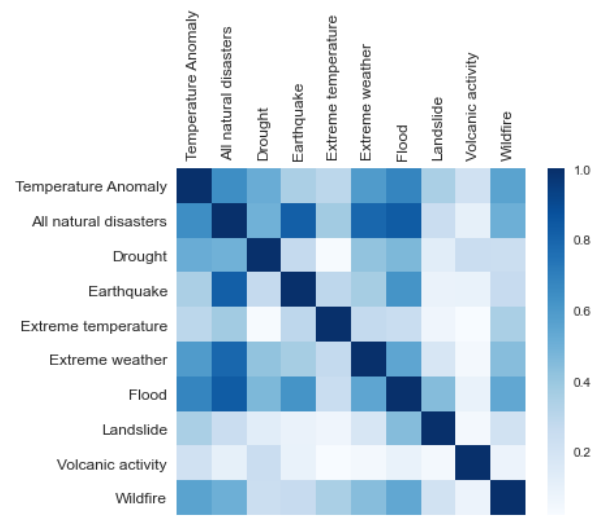

Fig. 16. Pearson Correlation Heatmap for dependencies between economic damages and temperature anomaly.

Kendall correlation values and heatmap showcase a different perspective, with a slightly lower positive correlation between the occurrences of disasters and temperature anomaly while maintaining a mixed (both positive and negative) correlation coefficient with respect to economic damages. Spearman correlation values and heatmap represent a stronger positive correlation than Pearson with major disaster types such as floods and earthquakes. Concerning economic damage, the correlation is more positive compared to other coefficients.

The custom-built dataset incorporating the analysis factors allowed for the creation of a visualization dashboard with interpretable charts, to provide valuable insights required for disaster planning. Initially, we had depicted basic graphs such as year-wise global distribution of deaths and country-wise share of losses. For a further distinguishing view, various types of disasters were included, majorly earthquakes, floods, droughts, volcanoes, landslides and wildfires. The damages caused by each disaster were assessed to rank them according to their impact. This resulted in the scope for performing regional analyses and grouping them according to socio-economic metrics. To shed some light on calamities causing gradual but immense repercussions such as food shortage, we emphasized making a comparative study of news coverage among various disaster types.

A quantitative analysis made using the charts assisted us to infer region and type-specific results for a better insight. Since 1900, Haiti has been the country with the highest death toll of 2,27,000 due to disasters mostly attributed to the intense Earthquakes followed by Indonesia with 1,85,000 deaths. There occurred 10,271 disasters worldwide between the years 1970 and 2010, with 23.7% and 16.9% of them occurring in America and Caribbean countries respectively and Asia being the worst-affected continent, accounting for roughly 40% of the total. Around 9 million people were domestically displaced in India in 2012 as a result of Himalayan flash floods, whereas the count stood at 19 million in China in 2008 as a result of the Great Sichuan earthquake. The Caribbean is the region where disaster damage has accounted for the highest percentage of GDP loss on average, topping 8% seven times. Central America comes in second, with disaster damage exceeding 8% of GDP twice.

While each disaster can be ranked through impact, it is observed that floods are the most common natural disaster by far, accounting for 43% of all recorded events. Predominantly, earthquakes, floods and droughts have resulted in a large number of deaths. Over the last few decades, high death tolls tend to result from large earthquake events, two such years being 2004 and 2010 where earthquakes accounted for 93% and 69% of total deaths respectively. Food shortages are responsible for most deaths and afflict higher number of casualties per incident occurred, but they commence in an inconspicuous manner as opposed to a sudden volcanic explosion or earthquake. As a result, food shortages are covered by news media only 3% of the time while a comparatively indulgent 30% of earthquakes and volcanic events get their time in the spotlight. To analyze the impact of each type of disaster and map them with a particular region, for instance, India, further quantitative analysis is performed and the results are interpreted. Subsequently, performing an exploratory data analysis on global temperature and disaster data indicated the suitability of the correlation approach from the available statistical methods. We studied the significance of each correlation coefficient and the varied results they produced with the dataset given. The interpretations made from correlation analysis reflected that climate change indeed enhances the disaster risk and occurrence, as well as economic losses, post the disaster.

## 5. Conclusion

The research aimed to explore dependencies between various disaster-related factors and use the interpretations for better awareness of natural calamities and resource planning during emergency response. The disaster data corpus was created and integrated with global temperatures and economic losses to develop a dashboard of visualizations and perform correlation analysis. We inferred that disasters affect those in poverty the most, high death tolls tend to be centered in low-to-middle income countries without the infrastructure to protect and respond to events. There is notable year-to-year variability in costs ranging from 0.15% to 0.5% of global GDP. In recent decades there has been no clear trend in damages when we take into account economic growth over this period. In India, the economic damages due to floods are the costliest, causing 63% of the total, followed by cyclones (19%), earthquakes (10%) and droughts (5%). In terms of human casualties, the earthquake is the most lethal disaster in India with 33% of casualties, followed by floods with 32%. The type of disaster matters to how newsworthy people find it to be. The news attention does not reflect the severity and number of people killed or affected by a natural disaster. Food shortages, for example, result in the most casualties and affect the most people per incident but their onset is more gradual than that of a volcanic explosion or sudden earthquake. Hence, nearly 40 thousand deaths due to food shortage have to be reported to receive the same newsworthiness as opposed to 2 deaths due to an earthquake.

Correlation analysis by medium variables was conducted using Pearson, Kendall, and Spearman statistical methods. Each medium variable and correlation was analyzed to be significant. Strong positive correlations were shown between the frequency of natural disasters and global warming and associated climate change scenarios, as environmentalists have been warning us. In terms of inferring the effect of disasters on the country's economy, a positive correlation was exhibited by Pearson and Spearman coefficients, while Kendall exhibited positive values for common and severely damaging types of disasters and varied values were observed for rare events. Therefore, developing countries are particularly vulnerable to catastrophes triggered by climate change and they need substantial economic support in the domain of disaster relief.

The future scope of the study includes the formulation of a disaster risk index by combining relevant factors from satellite imagery and social media data specific to each region to identify vulnerable regions using data patterns. Therefore, we can imply that the meaningful analysis of various types of disaster data pre and post of occurrence lays the foundation for disaster relief.

**References**

[1]  *EM-DAT (Emergency Events Database).*

[2]  *Earth Surface Temperature Data.*

[3]  *Recommendations on the Role of Official Statistics in Measuring Hazardous Events and Disasters.*

[4]  G. Lizarralde, L. Bornstein, M. Robertson, K. Gould, B. Herazo, A.-M. Petter, H. Páez, J. H. Díaz, A. Olivera and G. González, "Does climate change cause disasters? How citizens, academics, and leaders explain climate-related risk and disasters in Latin America and the Caribbean," *International Journal of Disaster Risk Reduction,* vol. 58, p. 102173, 2021.

[5]  U. N. O. f. D. R. R. (. Centre for Research on the Epidemiology of Disasters, *The Human Cost of Disasters. An Overview of the Last 20 Years (2000–2019),* CRED Geneva, Switzerland, 2020.

[6]  *Weather-related disaster statistics.*

[7]  M. Coronese, F. Lamperti, K. Keller, F. Chiaromonte and A. Roventini, "Evidence for sharp increase in the economic damages of extreme natural disasters," *Proceedings of the National Academy of Sciences,* vol. 116, p. 21450–21455, 2019.

[8]  T. Geiger and A. Stomper, "Rising economic damages of natural disasters: Trends in event intensity or capital intensity?," *Proceedings of the National Academy of Sciences,* vol. 117, p. 6312–6313, 2020.

[9]  *Disaster Risk and Vulnerability: The Role and Impact of Population and Society.*

[10] *Greater Impact: How Disasters Affect People of Low Socioeconomic Status.*

[11] J. Hendler, "Data integration for heterogeneous datasets," *Big data,* vol. 2, p. 205–215, 2014.

[12] B. Dageville, T. Cruanes, M. Zukowski, V. Antonov, A. Avanes, J. Bock, J. Claybaugh, D. Engovatov, M. Hentschel and J. Huang, "The snowflake elastic data warehouse," in *Proceedings of the 2016 International Conference on Management of Data*, 2016.

[13] P. Michele, F. Fallucchi and E. W. De Luca, "Create dashboards and data story with the Data & Analytics frameworks," in *Research Conference on Metadata and Semantics Research*, 2019.

[14] E. Bisong, "Matplotlib and Seaborn," in *Building machine learning and deep learning models on google cloud platform*, Springer, 2019, p. 151–165.

[15] W. Han, C. Liang, B. Jiang, W. Ma and Y. Zhang, "Major natural disasters in China, 1985–2014: occurrence and damages," *International journal of environmental research and public health,* vol. 13, p. 1118, 2016.

[16] S. Banholzer, J. Kossin and S. Donner, "The impact of climate change on natural disasters," in *Reducing disaster: Early warning systems for climate change*, Springer, 2014, p. 21–49.

[17] M. Kuncoro and J. Sartohadi, "Physical characteristics and disaster risk perception correlation at Bantul regency," *Economic Journal of Emerging Markets,* vol. 4, p. 76–86, 2012.

[18] V. Thomas, J. R. G. Albert and C. Hepburn, "Contributors to the frequency of intense climate disasters in Asia-Pacific countries," *Climatic Change,* vol. 126, p. 381–398, 2014.

[19] M. Berlemann and M. F. Steinhardt, "Climate change, natural disasters, and migration—a survey of the empirical evidence," *CESifo Economic Studies,* vol. 63, p. 353–385, 2017.

[20] V. F. Pisarenko and M. V. Rodkin, Statistical analysis of natural disasters and related losses, Springer, 2014.

**Authors' Profiles**

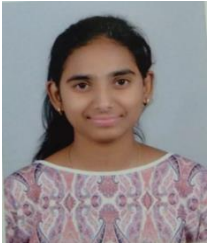

**Sukeerthi Mandyam,** is an Undergraduate student in the Department of Computer Science and Engineering at Sri Sivasubramaniya Nadar College of Engineering, Chennai. She has developed projects based on data structures, database systems, web development, object-oriented design and real-life ETL cycle using data warehousing and automation techniques. Her research interests include data analytics, data warehousing, machine learning, and image processing.

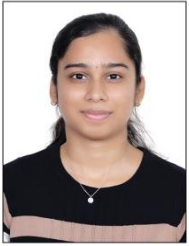

**Shanmuga Priya M.G**, is an Undergraduate Student in the Department of Computer Science and Engineering at Sri Sivasubramaniya Nadar College of Engineering, Chennai. She has worked on projects for managing databases by creating and connecting through Rest APIS, and creating user interfaces using React JS. She explores problem solving in the field of artificial intelligence, NLP, and Deep Learning.

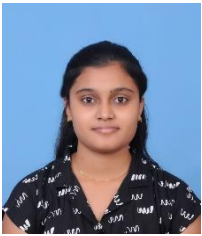

**Shalini S,** is an Undergraduate Student in the Department of Computer Science and Engineering at Sri Sivasubramaniya Nadar College of Engineering, Chennai. She has worked on real-time projects that involve latest technologies such as Ruby, Vue.js, Redis and MySQL. Her research interests include computer vision, image processing, application development and backend processing techniques.

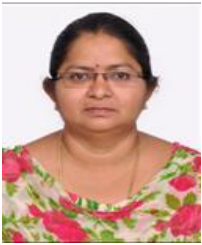

**Kavitha S**, is an Associate Professor in the Department of Computer Science and Engineering at Sri Sivasubramaniya Nadar College of Engineering, Chennai. She has 23 years of teaching experience, including 15 years of research experience in diverse domains namely bioinformatics, data mining, machine learning, soft computing and image analysis in digital and medical images, deep learning etc. She is a member of IEEE and the Machine Learning Research Group of SSN, a life member of Computer Society of India (CSI), and Indian Society for Technical Education (ISTE).  She has more than 70 publications in reputed Journals, National and International Conferences.